\documentclass[11pt]{article}
\usepackage[english]{babel}
\usepackage[utf8]{inputenc}
\usepackage[T1]{fontenc}
\usepackage{amsmath}
\usepackage{amssymb}
\usepackage{amsfonts}
\usepackage{amsthm}
\usepackage{bm}
\usepackage{hyperref}
\usepackage{upgreek}

\newcommand{\N}{\mathbb{N}}

\renewcommand{\H}{\mathbb{H}}
\newcommand{\R}{\mathbb{R}}

\newcommand{\boC}{\mathcal{C}}

\newcommand{\boH}{\mathcal{H}}

\newcommand{\boL}{\mathcal{L}}

\newcommand{\boN}{\mathcal{N}}
\newcommand{\boO}{\mathcal{O}}
\newcommand{\boP}{\mathcal{P}}

\newcommand{\boU}{\mathcal{U}}

\newtheorem{lem}{Lemma}
\newtheorem{prop}{Proposition}

\newtheorem{thm}{Theorem}
\theoremstyle{definition}
\newtheorem*{merci}{Acknowledgments}

\theoremstyle{remark}

\begin{document}

\title{Smooth travelling-wave solutions to the inviscid surface quasi-geostrophic equation}
\author{
\renewcommand{\thefootnote}{\arabic{footnote}}
Philippe Gravejat\footnotemark[1]~ and Didier Smets\footnotemark[2]}
\footnotetext[1]{Universit\'e de Cergy-Pontoise, Laboratoire de Math\'ematiques (UMR 8088), F-95302 Cergy-Pontoise Cedex, France. E-mail: {\tt philippe.gravejat@u-cergy.fr}}
\footnotetext[2]{Laboratoire Jacques-Louis Lions, Universit\'e Pierre et Marie Curie, Bo\^ite Courrier 187, 75252 Paris Cedex 05, France. E-mail: {\tt smets@ann.jussieu.fr}}
\maketitle

\begin{abstract}
We construct families of smooth travelling-wave solutions to the inviscid surface quasi-geostrophic equation~\eqref{eq:sqg}. These solutions can be viewed as the equivalents for this equation of the vortex anti-vortex pairs in the context of the incompressible Euler equation. Our argument relies on the stream function formulation and eventually amounts to solving a fractional nonlinear elliptic equation by variational methods.
\end{abstract}

\section{Introduction}

We consider the inviscid surface quasi-geostrophic equation 
\begin{equation}
\label{eq:sqg}
\tag{SQG}
\begin{cases}
\partial_t \theta + u \cdot \nabla \theta = 0,\\
u = R^\perp \theta,
\end{cases}
\end{equation}
where $R$ is the Riesz transform, $\theta : \R^2 \times \R \to \R$ is called the active scalar and $u : \R^2 \times \R \to \R^2$ is the velocity field induced by $\theta.$ Since $u$ is divergence free, it is convenient to relate $u$ and $\theta$ through a stream function $\psi : \R^2 \times \R \to \R$ by the equations
$$
\begin{cases}
u = \nabla^\perp \psi,\\
(- \Delta)^\frac{1}{2} \psi = \theta.
\end{cases}
$$
The inviscid surface quasi-geostrophic equation first appeared as a limit model in the context of geophysical flows. It has been widely investigated since the seminal work~\cite{ConMaTa1} of Constantin, Majda and Tabak, which pointed out its formal mathematical analogies with the three dimensional Euler equation. The Cauchy problem for~\eqref{eq:sqg} is known to be extremely delicate, and large classes of initial data are expected to produce finite time singularities. Besides radially symmetric solutions, which are all stationary, the only examples of global smooth solutions we are aware of were recently provided by Castro, C\'ordoba and G\'omez-Serrano~\cite{CasCoGS1}. We also refer to~\cite{CasCoGS1} for an extensive bibliography on the Cauchy problem for~\eqref{eq:sqg}. Our main goal in this note is to provide an alternative construction of smooth families of global special solutions.

We focus on travelling-wave solutions to~\eqref{eq:sqg}. Up to a rotation, we may assume, without loss of generality, that these waves have a positive speed $c$ in the vertical direction $z$, so that 
$$
\theta(r, z , t) = \Theta(r, z - c t), \quad u(r, z, t) = U(r, z - c t), \quad \psi(r, z, t) = \Psi(r, z - c t),
$$
for some profile functions $\Theta$, $U$ and $\Psi$ defined on $\R^2$. In this setting, equation~\eqref{eq:sqg} may be recast as the orthogonality condition
\begin{equation}
\label{eq:orth}
\big\langle \nabla^\perp \Psi - c e_z , \nabla \Theta \big\rangle_{\R^2} = 0,
\end{equation}
with $e_z = (0, 1)$. In the context of the Euler equation, Arnold~\cite{Arnold1} remarked that any function of the form
$$
\Theta(r, z) = f(\Psi(r, z) - c r) 
$$
automatically satisfies the orthogonality condition~\eqref{eq:orth}, at least formally, so that the travelling-wave problem reduces to a nonlinear elliptic equation. In our context, the same idea would lead to the fractional equation
$$
(- \Delta)^\frac{1}{2} \Psi = f(\Psi - c r). 
$$
We study a slight variation of this idea, in particular in order to get away from the radially symmetric situation.  We first assume a mirror symmetry with respect to the $z$-axis, namely
\begin{equation}
\label{eq:symmetry}
\Psi(r, z) = - \Psi(- r, z),
\end{equation}
so that $U(r, z) = U(- r, z)$ and $\Theta(r, z) = - \Theta(- r, z)$. We next impose the ansatz 
\begin{equation}
\label{eq:arnoldbis}
\Theta(r,z) = \begin{cases} \phantom{-} f \big( \Psi(r, z) - c r - k \big) & \text{ if } r \geq 0,\\
- f \big( - \Psi(r, z) + c r - k) & \text{ if } r \leq 0, \end{cases} 
\end{equation}
where $f$ is a smooth profile, and $k$ a positive number, to be specified later. In order to avoid any ambiguity or singularity for $r = 0$, we shall impose that $f(s)$ vanishes whenever $s \leq 0$. The condition~\eqref{eq:arnoldbis} also enforces the orthogonality condition~\eqref{eq:orth}, and leads likewise to the equation
\begin{equation}
\label{eq:main}
(- \Delta)^\frac{1}{2} \Psi = \begin{cases} \phantom{-} f \big( \Psi(r, z) - c r - k \big) & \text{ if } r \geq 0,\\
- f \big( - \Psi(r, z) + c r - k) & \text{ if } r \leq 0. \end{cases}
\end{equation}
Equation~\eqref{eq:main} is variational. Under our previous assumptions, its solutions are critical points of the functional
$$
E(\Psi) := \frac{1}{2} \int_{\R^2} \Psi (- \Delta)^\frac{1}{2} \Psi - \int_{\H} F(\Psi - c r - k) - \int_{\R^2 \setminus \H} F(- \Psi + c r - k),
$$
where we have set $\H := \{ (r, z) \in \R^2 \text{ s.t. } r \geq 0 \}$ and $F(s) := \int_0^s f(x) \, dx$. We construct critical points of $E$ as minimizers on the so-called Nehari manifold. For that purpose, we now make precise our functional framework. We assume~\footnote{Note that these assumptions together imply that $\nu \geq 2$.}
\begin{align}
\label{eq:H1}
\tag{$H_1$}
& f \in \boC^\infty(\R, \R) \setminus \{ 0 \},\\
\label{eq:H2}
\tag{$H_2$}
& f(s) = 0, \quad \forall\ s \leq 0,\\
\label{eq:H3}
\tag{$H_3$}
& f'''(s) \geq 0, \quad \forall\ s \geq 0,\\
\label{eq:H4}
\tag{$H_4$}
& \exists \nu < 3 \text{ and } C > 0 \text{ s.t. } f(s) \leq C s^\nu, \quad \forall\ s \geq 0.
\end{align}
A typical example verifying these assumptions is given by any function $f$ with $f(0)=f'(0)=f''(0)=0$ and $f'''=g$ where $g\: (\neq 0)$ is smooth, non-negative, and compactly supported in $\R_+.$ 
Under these assumptions, the functional $E$ is well-defined and differentiable on the Hilbert space
$$
X := \big\{ \Psi \in L^4(\R^2) \text{ s.t. } \Psi \in \dot{H}^\frac{1}{2}(\R^2) \big\},
$$
endowed with the scalar product
\begin{equation}
\label{eq:scalarprod}
\big\langle \Psi_1, \Psi_2 \big\rangle_X := \int_{\R^2} \int_{\R^2} \frac{\big( \Psi_1(x) - \Psi_1(y) \big) \, \big( \Psi_2(x) - \Psi_2(y) \big)}{|x - y|^3} \, dx \, dy.
\end{equation}
The energy $E$ is invariant under the symmetry group generated by~\eqref{eq:symmetry}. It follows from the Palais principle of symmetric criticality \cite{Palais1} that any critical point of the restriction of $E$ to the space $X_{\rm sym}$ of invariant functions is also a critical point of $E$ on the entire space $X$. In the sequel, we therefore restrict our analysis to the space $X_{\rm sym}$. In that space, the energy $E$ reduces to the expression
$$
E(\Psi) = \frac{1}{2} \big\| \Psi \big\|_X^2 - 2 V(\Psi) := \frac{1}{2} \int_{\R^2} \Psi (- \Delta)^\frac{1}{2} \Psi - 2 \int_{\H} F(\Psi - c r - k).
$$
The Nehari manifold associated to $E$ is defined by
$$
\boN = \big\{ \Psi \in X_{\rm sym} \setminus \{ 0 \} \text{ s.t. } E'(\Psi)(\Psi) = 0 \big\},
$$
so that $\Psi \in \boN$ if and only if
$$
\int_{\R^2} \Psi (- \Delta)^\frac{1}{2} \Psi - 2 \int_{\H} f(\Psi - c r - k) \, \Psi = 0.
$$
We shall prove that the set $\boN$ is a non-empty $\boC^1$-submanifold of $X_{\rm sym}$ without boundary. Our main result is then

\begin{thm}
\label{thm:main}
Let $c$ and $k$ be two positive numbers, and $f$ be an arbitrary profile verifying the assumptions~\eqref{eq:H1}-\eqref{eq:H4}. The functional $E$ possesses a minimizer $\Psi$ on $\boN$. As a consequence, there exists a non-trivial smooth travelling-wave solution $\theta$ to~\eqref{eq:sqg} given by
$$
\theta(r, z, t) = \Theta(r, z - c t) = f \big( \Psi(r, z - c t) - c r - k \big),
$$
for all $(r, z) \in \H$, and which satisfies the symmetry
$$
\Theta(r, z) = - \Theta(- r, z) = \Theta(r, - z),
$$
for all $(r, z) \in \R^2$. The restriction of $\Theta$ to $\H$ is non-negative with compact support, and is decreasing with respect to $|z|$.
\end{thm}

In the context of the two-dimensional and axisymmetric three-dimensional Euler equations, related constructions were first carried out by Berger and Fraenkel~\cite{BergFra1} and Norbury~\cite{Norbury1}. Contrary to these works, we do not know whether the support restricted to $\boH$ of the profile $\Theta$ in Theorem~\ref{thm:main} is connected.

\section{Strategy of the proof}

We consider the minimization problem
\begin{equation}
\label{eq:probmini}
\tag{$\boP$}
\alpha := \inf \big\{ E(\Psi), \ \Psi \in \boN \big\}.
\end{equation}
For $\Psi \in X_{\rm sym}$, we denote by $\Psi^\dagger$ the unique function, which is equal to the positive part $\Psi_+$ of $\Psi$ within $\H$, and which belongs to $X_{\rm sym}$. Since the nonlinearity $f$ identically vanishes on the negative axis, a function $\Psi$ cannot belong to $\boN$ if $\Psi^\dagger \equiv 0$. On the other hand, we have

\begin{prop}
\label{prop:nehari}
The Nehari constraint $\boN$ is a non-empty $\boC^1$-submanifold of $X_{\rm sym}$. For any $\Psi \in X_{\rm sym}$ with $\Psi^\dagger \neq 0$, there exists a unique positive number $t_\Psi$ such that $t_\Psi \Psi \in \boN$. The value of $t_\Psi$ is characterized by the identity
\begin{equation}
\label{eq:col}
E(t_\Psi \Psi) = \max \big\{ E(t \Psi), \ t > 0 \big\},
\end{equation}
and any critical point of $E$ on $\boN$ is a non-trivial smooth solution to~\eqref{eq:main}. Moreover, we have
\begin{equation}
\label{eq:parenbas}
\beta := \inf \big\{ \|\Psi\|_X^2, \ \Psi \in \boN \big\} > 0,
\end{equation}
and for any $\Psi \in \boN$,
\begin{equation}
\label{eq:bound} 
\| \Psi \|_X^2 \leq 6 E(\Psi).
\end{equation}
In particular, the minimal value $\alpha \geq \beta/6$ is positive, and any minimizing sequence for $E$ on $\boN$ is bounded.
\end{prop}

We notice that $E(\Psi^\dagger) \leq E(\Psi)$. A related observation is

\begin{lem}
\label{lem:pos}
For any $\Psi \in \boN$, we have
$$
E \big( t_{\Psi^\dagger} \Psi^\dagger \big) \leq E(\Psi),
$$
the inequality being strict whenever $\Psi$ is not equal to $\Psi^\dagger$.
\end{lem}

We denote by $X_{\rm sym}^\dagger$ and $\boN^\dagger$ the subsets of functions $\Psi$ in $X_{\rm sym}$, respectively $\boN$, which satisfy $\Psi = \Psi^\dagger$. From Lemma~\ref{lem:pos}, we deduce that
$$
\alpha := \inf \big\{ E(\Psi), \ \Psi \in \boN^\dagger \big\},
$$
and we therefore restrict our attention in the sequel to the functions $\Psi$ in $X_{\rm sym}^\dagger$. 

For $\Psi \in X_{\rm sym}^\dagger$, we denote by $\Psi^\sharp$ its Steiner symmetrization with respect to the vertical variable $z$. We observe that $E(\Psi^\sharp) \leq E(\Psi)$. Similarly to Lemma~\ref{lem:pos}, we have

\begin{lem}
\label{lem:steiner}
For any $\Psi \in \boN^\dagger$, we have
$$
E \big( t_{\Psi^\sharp} \Psi^\sharp \big) \leq E(\Psi).
$$
\end{lem}

In view of the information gathered so far, we may restrict our attention to a minimization sequence $(\Psi_n)_{n \in \N}$ for~\eqref{eq:probmini} such that $\Psi_n = \Psi_n^\dagger = \Psi_n^\sharp$. By Proposition~\ref{prop:nehari}, this sequence is bounded. We claim

\begin{lem}
\label{lem:comp}
Let $c$ and $k$ be positive numbers. The mapping
$$
\Psi \mapsto \big( \Psi^\sharp - c r - k \big)^\dagger
$$
is compact from $X_{\rm sym}^\dagger$ into $L^p(\R^2)$ for any number $1 \leq p < 4$.
\end{lem}

Passing to a subsequence if necessary, we assume that $\Psi_n \rightharpoonup \Psi_*$ weakly in $X$, and $(\Psi_n - c r - k)^\dagger \to (\Psi_* - c r - k)^\dagger$ strongly in $L^p(\R^2)$, as $n \to +\infty$, for any $1 \leq p < 4$.

\begin{prop}
\label{prop:strongconv}
The convergence of $\Psi_n$ towards $\Psi_*$ is strong in $X$. In particular, $\Psi_*$ is a solution to the minimization problem~\eqref{eq:probmini}.
\end{prop}

We finally define $\Theta_*$ from $\Psi_*$ according to~\eqref{eq:arnoldbis}, and we complete the proof of Theorem~\ref{thm:main} by

\begin{prop}
\label{prop:compactsupport}
The function $\Psi_*$ is smooth on $\R^2$, and there exists a positive number $C$ such that
$$
\big| \Psi_*(r,z) \big| \leq \frac{C}{1 + |r| + |z|}, \quad \forall (r,z) \in \R^2.
$$
In particular, the function $\Theta_*$ has compact support in $\R^2 \setminus \{ (r, z) \in \R^2 \text{ s.t. } r = 0 \}$.
\end{prop}

\section{Details of the proofs}

\subsection{Proof of Proposition~\ref{prop:nehari}}

Les us fix $\Psi \in X_{\rm sym}$, with $\Psi^\dagger \neq 0$. Given a positive number $t$, we let
$$
g(t) := \frac{E'(t \Psi)(t \Psi)}{t^2} = \frac{1}{2} \big\| \Psi \big\|_X^2 - \frac{2}{t} \int_{\H} f \big( t \Psi^\dagger - c r - k \big) \, \Psi^\dagger.
$$
We claim that the map $t \mapsto g(t)$ has one and only one zero in $\R_+$. As a consequence of our assumptions on $f'''$, we first observe that 
\begin{equation}
\label{eq:dertierce}
z f'(z) - 2 f(z) \geq 0, \quad \text{and} \quad f(z) \geq f(z_0) \Big( \frac{z}{z_0} \Big)^2,
\end{equation}
for all $z,z_0$ such that $z \geq z_0 > 0$. Since $\Psi^\dagger \neq 0$, we infer that 
$$
g(t) \to - \infty,
$$
as $t \to + \infty$. On the other hand, there exists a positive number $K$ such that $f(z) \leq K z^3$, when $z \geq 0$. Hence, we have
\begin{equation}
\label{eq:sanremo}
\frac{1}{t} f(t \Psi^\dagger - c r - k) \, \Psi^\dagger \leq \frac{1}{t} f(t \Psi^\dagger) \, \Psi^\dagger \leq K t^2 \big( \Psi^\dagger \big)^4,
\end{equation}
and therefore,
$$
\lim_{t \to 0} g(t) = \frac{1}{2} \| \Psi \|_X^2 > 0.
$$
By continuity, there exists at least a positive number $t_\Psi$ such that $g(t_\Psi) = 0$. We claim that $g'(t_\Psi) < 0$, which ensures the uniqueness of $t_\Psi$. For that purpose, we compute 
\begin{equation}
\label{eq:roubaix}
\begin{split}
t^2 g'(t) & = 2 \int_\H \Big( f(t \Psi^\dagger - c r - k) - t \Psi^\dagger f'(t \Psi^\dagger - c r - k) \Big) \, \Psi^\dagger\\
& \leq - 2 \int_\H \Big( f(t \Psi^\dagger - c r - k) + (c r + k) f'(t \Psi^\dagger - c r - k) \Big) \, \Psi^\dagger,
\end{split}
\end{equation}
where the last inequality follows from~\eqref{eq:dertierce}. Since $g(t_\Psi)=0$, we obtain
$$
\int_\H f(t_\Psi \Psi^\dagger - c r - k) \, \Psi^\dagger = \frac{t_\Psi}{4} \big\| \Psi \big\|_X^2 > 0.
$$
The uniqueness of $t_\Psi$ results from the non-negativeness of $c$, $k$ and $f'$. The characterization~\eqref{eq:col} is then a consequence of the identity $t g(t) = \frac{d}{dt} E(t \Psi)$. 

For $\Psi \in \boN$, we next write
$$
E(\Psi) = E(\Psi) - \frac{1}{3} E'(\Psi)(\Psi) = \frac{1}{6} \big\| \Psi \big\|_X^2 + \frac{2}{3} \int_\H \Big( f(\Psi - c r - k) \Psi - 3 F(\Psi - c r - k) \Big).
$$
By integration of~\eqref{eq:dertierce}, we know that $z f(z) - 3 F(z) \geq 0$ when $z \geq 0$, which gives~\eqref{eq:bound}. In view of~\eqref{eq:sanremo}, and the fact that $\Psi \in \boN$, we also have
$$
\frac{1}{2} \big\| \Psi \big\|_X^2 = 2 \int_{\H} f(\Psi^\dagger - c r - k) \, \Psi^\dagger \leq 2 K \int_{\H} \big( \Psi^\dagger \big)^4 \leq C \big\| \Psi \big\|_X^4,
$$
where we have used the Sobolev embedding theorem. This yields~\eqref{eq:parenbas}, with $\beta := 1/(2 C)$. The 
positivity of $\alpha$ follows combining~\eqref{eq:parenbas} and~\eqref{eq:bound}.

The smoothness of $\boN$ is then a consequence of the implicit function theorem applied to the smooth mapping $\Xi : (t, \Psi) \mapsto E'(t\Psi)(\Psi)$, which is defined on the open set $\R_+^* \times \{ \Psi \in X_{\rm sym} \text{ s.t. } \Psi^\dagger \neq 0 \}$. Indeed, whenever $\Psi \in \boN$, we deduce as in~\eqref{eq:roubaix} that 
$$
\partial_t \Xi(1, \Psi) = 2 \int_\H \Big( f(\Psi^\dagger - c r - k) - \Psi^\dagger f'(\Psi^\dagger - c r - k) \Big) \, \Psi^\dagger < 0.
$$
Finally, any minimizer of $E$ on $\boN$ is a global minimizer of the function $\Psi \mapsto E(t_\Psi\Psi)$ 
on the open set $\{ \Psi \in X_{\rm sym} \text{ s.t. } \Psi^\dagger \neq 0 \}$. Therefore, using the definition of the Nehari manifold and the fact that $t_\Psi = 1$ for $\Psi \in \boN$, we conclude that
$$
E'(\Psi)(h) = E'(t_\Psi \Psi)(t_\Psi'(h) \Psi + t_\Psi h) = 0, 
$$
for all $h \in X_{\rm sym}$. \qed

\subsection{Proof of Lemma~\ref{lem:pos}}

Let us first remark that $\Psi^\dagger \neq 0$, when $\Psi \in \boN$. In view of~\eqref{eq:col}, and the fact that $\Psi \in \boN$, we know that
$$
E(\Psi) \geq E(t_{\Psi^\dagger} \Psi).
$$
On the other hand, since $F$ vanishes on the negative axis, it holds
$$
V(t_{\Psi^\dagger} \Psi) = V(t_{\Psi^\dagger} \Psi^\dagger). 
$$
Finally, we deduce from the definition~\eqref{eq:scalarprod} of the scalar product in $X$, and from the fact that $\Psi$ and $\Psi^\dagger$ coincide on the support of $\Psi^\dagger$, that
$$
\| t_{\Psi^\dagger} \Psi\|_X^2 \geq \| t_{\Psi^\dagger} \Psi^\dagger \|_X^2, 
$$
the inequality being strict whenever $\Psi \neq \Psi^\dagger$. The conclusion follows combining the previous three arguments. \qed

\subsection{Proof of Lemma~\ref{lem:steiner}}

Arguing exactly as in the proof of Lemma~\ref{lem:pos}, it suffices to establish that
$$
E(\Psi^\sharp) \leq E(\Psi),
$$
when $\Psi \in X_{\rm sym}^\dagger$.
Since the Steiner symmetrization only involves rearrangements of super-level sets, we first have 
$$
V(\Psi^\sharp) = V(\Psi).
$$
On the other hand, we claim that
\begin{equation}
\label{eq:normdec}
\| \Psi^\sharp \|_X^2 \leq \| \Psi \|_X^2.
\end{equation}
This was proved e.g. by Almgren and Lieb~\cite[Theorem 9.2]{AlgrLie1}. For the sake of completeness, we present below a related short proof.

We first observe that
$$
\| \Psi \|_X^2 = \lim_{t \to 0} \int_{\R^2} \int_{\R^2} \frac{|\Psi(x) - \Psi(y)|^2}{(|x - y|^2 + t^2)^\frac{3}{2}} \, dx \, dy.
$$
For a compactly supported function $\Psi$ and for a fixed positive number $t$, symmetrizing the last expression in $x$ and $y$ yields the identity
$$
\int_{\R^2} \int_{\R^2} \frac{|\Psi(x) - \Psi(y)|^2}{(|x - y|^2 + t^2)^\frac{3}{2}} = 2 \int_{\R^2} \int_{\R^2} \frac{\Psi(x)^2}{(|x - y|^2 + t^2)^\frac{3}{2}} \, dx \, dy - 2 \int_{\R^2} \int_{\R^2} \frac{\Psi(x) \Psi(y)}{(|x - y|^2 + t^2)^\frac{3}{2}} \, dx \, dy. 
$$
In the right-hand side above, the first integral is invariant by any rearrangement, since it only depends on the super-level sets of $\Psi$. The second integral is decreased by the Steiner symmetrization by virtue of the Riesz rearrangement inequality. Passing to the limit $t \to 0$ and using the density of compactly supported functions in $X$ yields the conclusion~\eqref{eq:normdec}. \qed

\subsection{Proof of Lemma~\ref{lem:comp}}

Let $T : \Psi \mapsto (\Psi - c r - k)^\dagger$. We first claim that $T$ maps $X_{\rm sym}^\dagger$ into itself. Since $T(\Psi) \in L^4(\R^2)$, we are reduced to prove that $T(\Psi) \in \dot{H}^\frac{1}{2}(\R^2)$. We introduce the set
$$
\Omega(\Psi) := \big\{ (r, z) \in \H \text{ s.t. } \Psi(r, z) > c r + k \big\} = \big\{ (r, z) \in \H \text{ s.t. } T(\Psi)(r, z) > 0 \big\}.
$$
In order to compute the double integral defining the $\dot{H}^\frac{1}{2}$-norm of $T(\Psi)$, we split $\H$ as $\Omega(\Psi) \cup \Omega(\Psi)^c$. For sake of simplicity, we write $\Omega$ instead of $\Omega(\Psi)$ in the sequel. By the Sobolev embedding theorem, we have
\begin{equation}
\label{eq:setborne}
\boL^2 \big( \Omega \big) \leq \frac{1}{k^4} \int_\H \Psi^4 \leq \frac{C}{k^4} \big\| \Psi \big\|_X^4.
\end{equation}
First, we check that
$$
\int_\Omega \int_\Omega \frac{|T(\Psi)(x) - T(\Psi)(y)|^2}{|x - y|^3} \, dx \, dy \leq 2 \int_\Omega \int_\Omega \frac{|\Psi(x) - \Psi(y)|^2 + c^2 |r(x) - r(y)|^2}{|x - y|^3} \, dx \, dy,
$$
and, using~\eqref{eq:setborne} and the Riesz rearrangement inequality,
$$
\int_\Omega \int_\Omega \frac{|r(x) - r(y)|^2}{|x - y|^3} \, dx \, dy \leq \int_\Omega \int_\Omega \frac{dx \, dy}{|x - y|} \leq 2 \pi^\frac{1}{2} \boL^2(\Omega)^\frac{3}{2} \leq C \big\| \Psi \big\|_X^6.
$$
Next, we have
$$
\int_{\Omega^c} \int_{\Omega^c} \frac{|T(\Psi)(x) - T(\Psi)(y)|^2}{|x - y|^3} \, dx \, dy = 0,
$$
and we write the last term as
\begin{equation}
\label{eq:poggio}
\int_\Omega \int_{\Omega^c} \frac{|T(\Psi)(x) - T(\Psi)(y)|^2}{|x - y|^3} \, dx \, dy = \int_\Omega |T(\Psi)(x)|^2 \int_{\Omega^c} \frac{dy}{|x - y|^3} \, dx.
\end{equation}
For each fixed $x \in \Omega$, let
$$
\boO_x := \Big\{ y \in \Omega^c \text{ s.t. } r(y) \geq r(x) + \Lambda(x) := r(x) + \frac{1}{2c} \big( \Psi(x) - c r(x) - k \big) \Big\}.
$$
We divide $\Omega^c$ as $\boO_x \cup (\Omega^c \setminus \boO_x)$. On the one hand, the definition of $\boO_x$ provides
\begin{equation}
\label{eq:alaphilippe}
\int_{\boO_x} \frac{dy}{|x - y|^3} \leq \frac{C}{\Lambda(x)}.
\end{equation}
On the other hand, for $x \in \Omega$ and $y \in \Omega^c \setminus \boO_x$, we can use the definition of $\boO_x$ in order to get
$$
0 \leq T(\Psi)(x) = \Psi(x) - c r(x) - k \leq \Psi(x) - \Psi(y) + c \big( r(y) - r(x) \big) \leq \Psi(x) - \Psi(y) + \frac{1}{2} T(\Psi)(x),
$$
so that 
\begin{equation}
\label{eq:sagan}
|T(\Psi)(x)| \leq 2 |\Psi(x) - \Psi(y)|.
\end{equation}
Combining~\eqref{eq:alaphilippe} and~\eqref{eq:sagan} in~\eqref{eq:poggio}, we deduce
$$
\int_\Omega \int_{\Omega^c} \frac{|T(\Psi)(x) - T(\Psi)(y)|^2}{|x - y|^3} \, dx \, dy \leq 4 \int_\Omega \int_{\Omega^c} \frac{|\Psi(x) - \Psi(y)|^2}{|x - y|^3} \, dx \, dy + C \int_\Omega |T(\Psi)(x)| \, dx,
$$
and we may additionally bound the last term in this sum as
$$
\int_\Omega |T(\Psi)(x)| \, dx \leq C \boL^2(\Omega)^\frac{3}{4} \big\| \Psi \big\|_{L^4} \leq C \big\| \Psi \big\|_X^4, 
$$
by invoking~\eqref{eq:setborne}, the Sobolev embedding theorem, and the fact that $T(\Psi) \leq \Psi$. Combining further all our estimates so far, we finally infer that
$$
\int_\Omega \int_{\Omega^c} \frac{|T(\Psi)(x) - T(\Psi)(y)|^2}{|x - y|^3} \, dx \, dy \leq C \big\| \Psi \big\|_X^2 \Big( 1 + \big\| \Psi \big\|_X^4 \Big),
$$
so that $T$ is well-defined from $X_{\rm sym}^\dagger$ into itself, and maps bounded sets into bounded sets.

We next turn to the compactness properties. Let $\Psi^\sharp$ be a Steiner symmetric function in $X_{\rm sym}^\dagger$. We claim that 
\begin{equation}
\label{eq:huy}
\int_{|z| \geq R} \big| T(\Psi^\sharp) \big|^2 \leq \frac{2^\frac{3}{2}}{\pi^\frac{3}{2} k^2 R} \big\| T(\Psi^\sharp) \big\|_X^2 \big\| \Psi^\sharp \big\|_{L^4}^2,
\end{equation}
for all $R \geq 1$. Indeed, let $\kappa$ be a positive number to be fixed later, and set $\boU_x := B(x, \kappa/R) \cap \Omega^c$ for all $x \in \H$. We estimate 
$$
\int_{|z(x)| \geq R} \big| T(\Psi^\sharp)(x) \big|^2 \, dx \leq \int_{|z(x)| \geq R} \frac{\kappa^3}{R^3 \boL^2(\boU_x)} \int_{\boU_x} \frac{|T(\Psi^\sharp)(x) - T(\Psi^\sharp)(y)|^2}{|x - y|^3} \, dy \, dx.
$$
As a consequence of the Steiner symmetry of $\Psi^\sharp$, and~\eqref{eq:setborne}, we have 
$$
\boL^1 \Big( \big\{ r \text{ s.t. } (r, z) \in \Omega^c \text{ for some } |z| \geq R \big\} \Big) \leq \frac{1}{2 R k^4} \big\| \Psi \big\|_{L^4}^4.
$$
We now fix $\kappa$ such that $\pi k^4 \kappa^2 = 2 \| \Psi \|_{L^4}^4$. As a consequence, computing the area of a disc minus a strip gives 
$$
\boL^2(\boU_x) \geq \frac{\pi \kappa^2}{R^2} - \frac{\| \Psi \|_{L^4}^4}{k^4 R^2} = \frac{\| \Psi \|_{L^4}^4}{k^4 R^2},
$$
so that~\eqref{eq:huy} follows.

Similarly, we have
\begin{equation}
\label{eq:redoute}
\int_{|r| \geq R} \big| T(\Psi^\sharp) \big|^2 \leq \frac{1}{(c R + k)^2} \big\| \Psi^\sharp \big\|_{L^4}^4.
\end{equation}
Indeed, we deduce from the definition of $\Omega$ that
$$
\int_{|r| \geq R} \big| T(\Psi^\sharp) \big|^2 \leq \boL^2 \big( \Omega \cap \{ (r, z) \in \H \text{ s.t. } |r| \geq R \} \big)^\frac{1}{2} \big\| \Psi^\sharp \big\|_{L^4}^2,
$$
and also 
$$
\boL^2 \big( \Omega \cap \{(r, z) \in \H \text{ s.t. } |r| \geq R \} \big) \leq \frac{1}{(c R + k)^4} \int_{\Omega \cap \{|r|\geq R\}} \big( \Psi^\sharp \big)^4.
$$
The conclusion then follows from the Rellich compactness theorem (at the local level) combined with the decay estimates in~\eqref{eq:huy} and~\eqref{eq:redoute}. \qed

\subsection{Proof of Proposition~\ref{prop:strongconv}}

It first follows from the definition of $\boN$ and Lemma~\ref{lem:comp} that
$$
\int_\H f(\Psi_* - c r - k) \, \Psi_* \geq \frac{\beta}{4}. 
$$
In particular, $(\Psi_* - c r - k)^\dagger \neq 0$, so that by Proposition~\ref{prop:nehari}, there exists a unique positive number $t_*$ such that $t_* \Psi_* \in \boN$. We shall prove that $t_* = 1$. Indeed, we have
$$\alpha = \lim_{n \to + \infty} E(\Psi_n) \geq \liminf_{n \to + \infty} E(t_*\Psi_n),$$
by Proposition~\ref{prop:nehari}, and 
$$\liminf_{n \to + \infty} E(t_* \Psi_n) \geq E(t_* \Psi_*) \geq \alpha,
$$
by Lemma~\ref{lem:comp}, and since $t_* \Psi_* \in \boN$. It follows that all these inequalities are equalities. In particular, we infer that $\lim_{n \to + \infty} \| \Psi_n \|_X^2 = \| \Psi_* \|_X^2$, from which the strong convergence of $\Psi_n$ towards $\Psi_*$ in $X$ follows. The latter implies that $\Psi_* \in \boN$, and therefore, that $t_* = 1$ and $E(\Psi_*) = \alpha$. \qed

\subsection{Proof of Proposition~\ref{prop:compactsupport}}

We already know that $T(\Psi_*) \in L^4(\R^2)$. Since the support of $T(\Psi_*)$ has finite measure, this implies that $T(\Psi_*)\in L^1(\R^2)$. Let $\Theta_*$ be defined through~\eqref{eq:arnoldbis}, where $\Psi$ is replaced by $\Psi_*$. It follows from~\eqref{eq:H4} that
$$
\Theta_* \in L^\frac{q}{\nu}(\R^2), \quad \forall\ \nu \leq q \leq 4.
$$
Consider next the function $\tilde{\Psi}_*$ given by the representation formula
$$
\tilde{\Psi}_* (x) := \frac{1}{2 \pi} \int_{\R^2} \frac{\Theta_*(y)}{|x - y|} \, dy.
$$
It follows from the weighted inequalities for singular integrals in~\cite[Chapter 5, Theorems 1 and 2]{Stein0} that
$$
\tilde{\Psi}_* \in \dot{W}^{1, q}(\R^2), \quad \forall 1 < q\leq \frac{4}{\nu}.
$$
Moreover, by the Hardy-Littlewood-Sobolev inequality, we have
$$
\tilde{\Psi}_* \in L^q(\R^2), \quad \forall 2 < q \leq \frac{4}{\nu - 2}.
$$
Hence, we deduce from standard interpolation that $\tilde{\Psi}_* \in X_{\rm sym}$. For $\varphi \in X_{\rm sym} \cap \boC_c^\infty(\R^2)$, a direct computation provides
$$
\big\langle \tilde{\Psi}_*, \varphi \big\rangle_X = \int_{\R^2} \Theta_* \varphi.
$$
On the other hand, since $\Psi_*$ is a critical point of $E$, we also have
$$
\big\langle \Psi_*, \varphi \big\rangle_X = \int_{\R^2} \Theta_* \varphi.
$$
By density of smooth compactly supported functions in $X_{\rm sym}$, it follows that $\tilde{\Psi}_* = \Psi_*$, and we may invoke a direct $L^p$-type bootstrap argument to deduce that $\Psi_*$ is bounded and uniformly continuous on $\R^2$. In order to deduce from a further bootstrap argument that $\Psi_*$ is smooth, we only need to check that the possible discontinuity at $r = 0$ introduced by the definition~\eqref{eq:arnoldbis} does not arise. This follows from the mirror symmetry assumption on $\Psi_*$, and the fact that $\Psi_*$ is already known to be uniformly continuous, so that $T(\Psi_*)$ vanishes in
an open strip containing the axis $r = 0$, and therefore $\Theta_*$ has at least the same regularity as $\Psi_*$.

It remains to compute the decay of $\Psi_*$. For that purpose, let $x \in \R^2$ be such that $|x| \geq 1$. We write
$$
2 \pi \Psi_*(x) = \int_{|x - y| \leq \frac{|x|}{2}} \frac{\Theta_*(y)}{|x-y|} \, dy + \int_{|x - y| \geq \frac{|x|}{2}} \frac{\Theta_*(y)}{|x-y|} \, dy.
$$
On the one hand, we estimate
\begin{equation}
\label{eq:rocheauxfaucons}
\int_{|x - y| \geq \frac{|x|}{2}} \frac{\Theta_*(y)}{|x - y|} \, dy \leq \frac{C}{|x|} \int_{\R^2} |T(\Psi_*)|^\nu \leq \frac{C}{|x|}.
\end{equation}
On the other hand, since $\Psi_*$ is Steiner symmetric, it follows from~\eqref{eq:huy} and~\eqref{eq:redoute} that
\begin{equation}
\label{eq:hautelevee}
\int_{|x - y| \leq \frac{|x|}{2}} T(\Psi_*)^2(y) \, dy \leq \frac{C}{|x|}.
\end{equation}
We infer from the Hölder inequality and~\eqref{eq:hautelevee} that
\begin{equation}
\label{eq:hans}
\int_{|x - y| \leq \frac{|x|}{2}} \frac{\Theta_*(y)}{|x - y|} \, dy \leq \bigg( \int_{|z| \leq \frac{|x|}{2}} \frac{dz}{|z|^\frac{5}{3}} \bigg)^\frac{3}{5} \, \bigg( \int_{|y - x| \leq \frac{|x|}{2}} T(\Psi_*(y))^\frac{5 \nu}{2} \, dy \bigg)^\frac{2}{5} \leq \frac{C}{|x|^\frac{1}{5}} \| \Psi_* \|_{L^\infty}^{\nu - \frac{4}{5}}.
\end{equation}
Combining~\eqref{eq:rocheauxfaucons} and~\eqref{eq:hans}, we deduce that
$$
|\Psi_*(x)| \leq \frac{C}{(1 + |x|)^\frac{1}{5}}, \quad \forall x \in \R^2.
$$
In view of the positive cut-off level $k$ entering in the definition of $T$, the latter implies that $T(\Psi_*)$, and therefore $\Theta_*$, have compact support. In turn, this implies that the left-hand side of~\eqref{eq:hans} vanishes for $|x|$ large. The conclusion follows from~\eqref{eq:rocheauxfaucons}. \qed

\begin{merci}
D.S. is partially supported by  grant ANR-14-CE25-0009-01 of the Agence Nationale de la Recherche. 
\end{merci} 
\bibliographystyle{plain}
\bibliography{Bibliogr}

\end{document}